
\magnification = \magstep 1

\baselineskip 15pt
\pageno=1

\def\hhrule#1#2{\kern-#1
   \hrule height#1 depth#2 \kern-#2 }
\def\hvrule#1#2{\kern-#1{\dimen0=#1
    \advance\dimen0 by#2\vrule width\dimen0}\kern-#2 }
\def\makeblankbox#1#2{\setbox0=\hbox{A}
\hbox{\lower\dp0\vbox{\hhrule{#1}{#2}%
   \kern -#1
   \hbox to \wd0{\hvrule{#1}{#2}%
     \raise\ht0\vbox to #1{}
     \lower\dp0\vtop to #1{}
     \hfil\hvrule{#2}{#1}}
  \kern-#1\hhrule{#2}{#1}}}}
\def\qed{\makeblankbox{0pt}{.3pt}}

\def\pra{\par}
\def\Hom{{\rm Hom}}
\def\Aut{{\rm Aut}}
\def\im{{\rm Im}}
\def\ker{{\rm Ker}}
\def\res{{\rm Res}}
\def\fp{{\bf F}_p}
\def\ftwo{{\bf F}_2}
\def\cat{{\cal C}}
\def\dat{{\cal D}}
\def\cats{{\cat_s}}
\def\catu{{\cat_u}}
\def\presc#1#2{{\vphantom{#2}^{#1}\negthinspace#2}}

\def\br{1}
\def\ce{2}
\def\wdd{3}
\def\ep{4}
\def\ev{5}
\def\se{6}
\def\smsw{7}

\font \smallfont=cmr8 at 8pt
\font \smallbold=cmbx8 at 8pt
\font \smallmath=cmmi8 at 8pt

\phantom{}
\vskip 1 truecm

\centerline{\bf ON UNIVERSALLY STABLE ELEMENTS}
\bigskip
\centerline{Ian J. Leary, Bj\"orn Schuster, and Nobuaki Yagita}
\bigskip
{{\narrower\smallskip\noindent{\bf Abstract.} 
{\smallfont 
We show that certain subrings of the cohomology of a finite
{\smallmath p}-group {\smallmath P} may be realised as the images of
restriction from suitable virtually free groups.  We deduce that the 
cohomology of {\smallmath P} is a finite module for any such subring.
Examples include the ring of \lq universally stable elements' defined
by Evens and Priddy, and rings of invariants such as the mod-2 
Dickson algebras.}
\smallskip}}

\bigskip

Let $P$ be a finite $p$-group, and let $\catu$ be the category whose
objects are the subgroups of $P$, with morphisms all injective group
homomorphisms.  Let $\cat$ be any subcategory of $\catu$ such that $P$
is an object of $\cat$, and such that for any object $Q$ of $\cat$, 
the inclusion of $Q$ in $P$ is a morphism in $\cat$.  Let
$H^*(\,\cdot\,)$ stand for mod-$p$ 
group cohomology, which may be viewed as a contravariant functor from 
$\catu$ to $\fp$-algebras.  We shall study the limit $I(P,\cat)$ of
this functor:  
$$I(P,\cat)= \lim_{Q\in \cat}H^*(Q).$$
Given our assumptions on $\cat$, we may identify $I(P,\cat)$ with a
subring of $H^*(P)$.  In the final remarks we discuss generalizations
of our results in which most of the conditions that we impose
upon $P$, $\cat$, and $H^*(\,\cdot\,)$ are weakened.  

The classical case of this construction occurs in Cartan-Eilenberg's
description of the image of the cohomology of a finite group in the
cohomology of its Sylow subgroup as the \lq stable elements\rq\ [\ce].  
Let $G$ be a finite group with 
$P$ as its Sylow subgroup, and let $\cat_G$ be the subcategory of 
$\catu$ containing all the objects, but with morphisms only those
homomorphisms $Q$ to $Q'$ induced by conjugation by some element of 
$G$.  Then the image, $\im(\res^G_P)$, of $H^*(G)$ in $H^*(P)$ is 
equal to $I(P,\cat_G)$.  

Rings of invariants also arise in this way.  If $\cat$ is a category
whose only object is $P$, with morphisms a subgroup $H$ of $\Aut(P)$,
then $I(P,\cat)$ is just the subring $H^*(P)^H$ of invariants under
the action of $H$.  

Another case considered already, which motivated our work,
is the ring of universally stable elements defined by Evens-Priddy
in [\ep].  Let $\cats$ be the subcategory of $\catu$ generated by all 
subcategories of the form $\cat_G$ as defined above.  Then 
$I(P,\cats)$ is the subring $I(P)$ of $H^*(P)$ introduced in [\ep],
consisting of those elements of $H^*(P)$ which are in the image of 
$\res^G_P$ for every finite group $G$ with Sylow subgroup $P$.  

A fourth case of interest is $I(P,\catu)$, which might be viewed as the
elements of $H^*(P)$ which are \lq even more stable' than the
elements of $I(P,\cats)$.  It is easy to see that in general $\cats$
is strictly contained in $\catu$.  For example, the endomorphism
monoid $\Hom_\cats(P,P)$ of $P$ is the subgroup of $\Aut(P)$ generated
by elements of order coprime to $p$, whereas $\Hom_\catu(P,P)$ is the
whole of $\Aut(P)$.  Our main result is the following theorem.   

\proclaim Theorem 1.  Let $P$ be a finite $p$-group, and let $\cat$ be
any subcategory of $\catu= \catu(P)$ satisfying the conditions
stated in the first paragraph.  Then there exists a discrete group $\Gamma$
containing $P$ as a subgroup such that: 
\item{a)} $\im(\res^\Gamma_P)$ is equal to $I(P,\cat)$;  
\pra 
\item{b)} $(\ker(\res^\Gamma_P))^2$ is trivial; 
\pra 
\item{c)} $\res^\Gamma_P$ induces an isomorphism from
$H^*(\Gamma)/\sqrt 0$ to $I(P,\cat)/\sqrt 0$;  
\item{d)} $\Gamma$ is virtually free.  More precisely, 
$\Gamma$ has a free normal subgroup of index dividing~$|P|!$.  

If $\Gamma'$ is a free normal subgroup of $\Gamma$ of finite index, 
then $P$ maps injectively to the finite group $\Gamma/\Gamma'$, so 
by the Evens-Venkov theorem [\ev], $H^*(P)$ is a finite module for
$H^*(\Gamma/\Gamma')$ and hence {\it a fortiori\/} for $H^*(\Gamma)$.
Thus one obtains the following corollary.  

\proclaim Corollary 2.  Let $P$ and $\cat$ be as in the statement of
Theorem~1.  Then $H^*(P)$ is a finite module for its  subring
$I(P,\cat)$.  

The case $\cat = \cats$ is Theorem~A of [\ep].  Our result is
stronger, since it applies to categories such as $\catu$ itself, and
our proof is more elementary.  There is an even shorter proof of
Corollary~2 however, which is to deduce it from the following 
simpler theorem.

\proclaim Theorem 3.  Let $P$ be a finite $p$-group, and let $G$ be
the symmetric group on a set $X$ bijective with $P$.  Regard $P$ as a
subgroup of $G$ via a Cayley embedding (or regular permutation
representation).  Then $\im(\res^G_P)$ is contained in $I(P,\catu)$.  

To deduce Corollary~2 from Theorem~3, note that for any $\cat$ as
above, one has 
$$\im(\res^G_P)\subseteq I(P,\catu)\subseteq I(P,\cat)\subseteq
H^*(P),$$
and $H^*(P)$ is a finite module for $\im(\res^G_P)$ by the
Evens-Venkov theorem.


\noindent
{\it Proof of Theorem 3.\/} We deduce Theorem~3 from the following
group-theoretic lemma.  

\proclaim Lemma 4.  Let $Q\leq P \leq G$ be as in the statement of
Theorem~3, and let $\phi$ be any injective homomorphism from $Q$ to
$P$.  Then there exists $g\in G$ such that for all $q\in Q$,  
$\phi(q) = gqg^{-1}$.  

\noindent
{\it Proof.\/} Fix a bijection between $P$ and the set $X$ permuted by
$G$.  This fixes an embedding $i_P$ of $P$ in $G$.  Let $i_Q$ be the 
induced inclusion of $Q$ in $G$.  Write $\presc{i_P}X$ for $X$ viewed
as a $P$-set.  Thus $\presc{i_P}X$ is a free $P$-set of rank one.  
There are two ways to view $X$ as a $Q$-set, either via 
$i_Q$ or $i_P\circ\phi$.  The $Q$-sets $\presc{i_Q}X$ and
$\presc{i_P\circ\phi}X$ are both free of rank equal to the index,
$|P:Q|$, of $Q$ in $P$.  Let $g$ be an isomorphism of $Q$-sets between
$\presc{i_Q}X$ and $\presc{i_P\circ\phi}X$.  Then $g$ is an element of $G$
having the required property, because for each $x\in X$ and $q\in Q$, 
$g\cdot q\cdot x = \phi(q)\cdot g\cdot x$.  \qed

Returning now to the proof of Theorem~3, any morphism in $\catu$ 
factors as the composite of an isomorphism followed by an inclusion.  
Thus it suffices to show that for $\phi$ as in Lemma~4, $\res^G_Q$ and
$\phi^*\circ\res^G_P$ are equal.  Writing $c_g$ for the automorphism
of $G$ given by conjugation by $g$, we have shown that there exists
$g$ such that $c_g\circ i_Q= i_P\circ\phi$.  But $c_g^*$ is the
identity map on $H^*(G)$, and hence $i_Q^*= \phi^*\circ i_P^*$ as
required.  \qed

This completes the proofs of all of our statements except for
Theorem~1.  For the proof of Theorem~1 we recall the following 
theorem (see for example [\wdd], I.7.4 or IV.1.6).

\proclaim Theorem 5.  Let $\Gamma$ be a group that acts simplicially
(i.e., without reversing any edges) on a tree with all stabilizer
groups of order dividing a fixed integer $M$.  Then there is a
homomorphism from $\Gamma$ to the symmetric group on $M$ letters whose
kernel, $K$, is torsion-free.  Since $K$ acts freely, simplicially on
the tree, it follows that $K$ is a free group.  

In fact, the short direct proof of Theorem~3  
is based on some of the ideas in the proof of Theorem~5 given in [\wdd].

\noindent
{\it Proof of Theorem 1.\/} The group $\Gamma$ will be constructed as
the fundamental group of a graph of groups (see [\wdd], I.3, [\se],
I.5, or [\br], VII.9 for the definitions and basic theorems).  
Let $Q_1,\ldots,Q_M$ be
the objects of $\cat$, and let $\phi_1,\ldots,\phi_N$ be the
morphisms of $\cat$.  
Define a function $m$ so that the domain of
$\phi_i$ is $Q_{m(i)}$.  Now let
$\Gamma$ be the group generated by the elements of $P$ and new
elements $t_1,\ldots,t_N$ subject to all relations that hold in $P$,
together with the relations 
$$t_iqt_i^{-1} = \phi_i(q),$$
for all $i\in \{1,\ldots,N\}$ and all $q \in Q_{m(i)}$.  
Thus $\Gamma$ is the fundamental group of a graph of groups with one
vertex and $N$ edges.  The vertex group is of course $P$ and the $i$th
edge group is $Q_{m(i)}$.  The two maps from the $i$th edge group to
the vertex group (corresponding to its initial and terminal ends) are
the inclusion and $\phi_i$.  

The group $\Gamma$ as defined above has the following properties (see
any of the references listed above):  
$P$ is a subgroup of $\Gamma$; for each $i$, the homomorphism 
$\phi_i\colon Q_{m(i)} \rightarrow P$ is inner in $\Gamma$ (i.e., is
induced by conjugation by the element $t_i$); $\Gamma$ acts
simplicially on a tree $T$, with one orbit of vertices and $N$ orbits of
edges, with $P$ being a vertex stabilizer and $Q_{m(i)}$ being the
stabilizer of some edge in the $i$th orbit.  The quotient $T/\Gamma$ 
is the graph used in defining $\Gamma$.  

Recall from [\br], VII.7--VII.9 that for any $\Gamma$-CW-complex $X$,
there is a spectral sequence, with $E_1^{p,q} = \bigoplus_\sigma
H^q(\Gamma_\sigma)$, where the sum is over a set of orbit
representatives of $p$-cells in $X$.  For coefficients in a ring with 
trivial $\Gamma$-action (such as the field of $p$ elements), this is a 
spectral sequence of rings.  When $X$ is acyclic the spectral sequence
converges to a filtration of $H^{p+q}(\Gamma)$.  
We apply this spectral sequence in the case when $X=T$.  In this 
case 
$$E_1^{0,q}\cong H^q(P),\quad E_1^{1,q}\cong \bigoplus_{i=1}^N
H^q(Q_{m(i)}),$$ 
and $E_1^{p,q}=0$ for $p>1$.  Under this isomorphism the differential 
$d_1:E_1^{0,q}\rightarrow E_1^{1,q}$ has $i$th coordinate
$\res^P_{Q_{m(i)}}-\phi^*_i$, and so $E_2^{0,*}$ is isomorphic to
$I(P,\cat)$.  The fact that $E_2^{p,q}=0$ for $p>1$ implies that the
spectral sequence collapses at the $E_2$-page.  The edge homomorphism
from $E_\infty^{0,*}$ to $H^*(P)$ may be identified with
$\res^\Gamma_P$ (consider the map of spectral sequences induced by the
inclusion of the vertex set of the tree in the whole tree, viewed as a
map of $\Gamma$-spaces), and so a) is proved.  For b), note that
since $E_2^{p,q}=0 $ for $p>1$, elements of $E_2^{1,*}$ uniquely
determine elements of $H^*(\Gamma)$, and the product of any two such
elements is zero in $H^*(\Gamma)$.  Since 
$\ker(\res^\Gamma_P)$ may be identified with $E_2^{1,*}$, 
b) follows, and c) follows immediately from b).  Finally, d) follows
from Theorem~5 stated above.  
\qed

\bigskip
\goodbreak
\noindent
{\bf Remarks. } 1) There are alternatives to using the equivariant
cohomology spectral sequence in the proof of Theorem~1, but following
a suggestion of the referee we decided to explain just one method in 
detail in the proof.  Since the spectral sequence has only two
non-zero rows it is essentially just a long exact sequence.  This long
exact sequence may be obtained by applying $H^*(\Gamma;\,\cdot\,)$ to 
the augmented chain complex for the tree $T$, modulo an application of
the Eckmann-Shapiro lemma.  We felt, however, that the ring
structure of $H^*(\Gamma)$ is more easily understood in terms of the
spectral sequence.  

\medskip\noindent
2) We believe that $I(P,\catu)$ has some
advantages over $I(P,\cats)$.  Both of these rings enjoy the 
finiteness property stated in Corollary~2.  To compute $I(P,\cats)$ 
one needs to know something about the $p$-local structure of all
groups with Sylow subgroup $P$, whereas $I(P,\catu)$ requires only
knowledge of $P$.  

\medskip\noindent 
3) On the other hand, $I(P,\catu)$ does not retain much information
concerning $P$.  Let $W(P)$ be the variety of all ring homomorphisms  
from $I(P,\catu)$ to an algebraically closed field $k$ of characteristic $p$.  
Then $W(P)$ is determined up to homeomorphism by the $p$-rank of $P$:
If $P$ has $p$-rank $n$, then $W(P)$ is homeomorphic to
$k^n/GL_n(\fp)$, and if $E$ is an elementary abelian subgroup of $P$ 
of rank $n$, then the induced map from $W(E)$ to $W(P)$ is an
homeomorphism.  These assertions concerning $W(P)$ follow easily from 
Quillen's theorem describing the variety of homomorphisms from $H^*(P)$ to
$k$ (see for example [\ev], chap.\ 9).  Note that this the only place 
where we use Quillen's theorem.  

\medskip\noindent
4) The definitions and theorems that we state remain valid if $P$ is
any finite group.  We restrict to the case when $P$ is a $p$-group 
only because this is the case occurring naturally in the work of 
Cartan-Eilenberg and Evens-Priddy.  

\medskip\noindent
5) The reader may have noticed that Theorems 1~and~3 work perfectly
well for cohomology with coefficients in any ring $R$ (viewed as a
trivial $P$-module).  
Corollary~2 is valid for cohomology with coefficients in any
ring $R$ for which the Evens-Venkov 
theorem holds (see [\ev], 7.4 for a general statement).  

\medskip\noindent
6) The easiest way to relax the restrictions on the category $\cat$ is
to consider arbitrary finite categories with objects finite groups and 
morphisms injective group homomorphisms (it is 
unhelpful to view the groups as subgroups of a single group if the
inclusion maps are not in the category).  Define $I(\cat)$ to be 
the limit over this category and for any group $\Gamma$, define 
$\dat(\Gamma)$ to be the category of finite subgroups of $\Gamma$, 
with morphisms inclusions and conjugation by elements of $\Gamma$.  
Then one obtains 
\proclaim Theorem $1'$.  Let $\cat$ be a connected finite category of
finite groups and injective homomorphisms.  Then there exists a
discrete group $\Gamma$ and a natural transformation from 
$\cat$ to $\dat(\Gamma)$ such that $\Gamma$ and the induced map 
from $H^*(\Gamma)$ to $I(\cat)$ satisfy properties a) to d) of 
Theorem~1.  

Recall that a category is said to be connected if the equivalence
relation on objects generated by \lq there is a morphism from $Q$ to
$Q'$' has exactly one class.  Note that there cannot be a direct
analogue of Theorem~1 unless the category $\cat$ is connected, since
the degree zero part of $I(\cat)$ is an $\fp$ vector space of
dimension the number of components of $\cat$, whereas
$H^0(\Gamma)\cong\fp$.  The proof of Theorem~$1'$ is very similar to
the proof of Theorem~1, except that one creates a graph of groups with
one vertex for every object of $\cat$.  The restriction to connected
categories is not serious, since given any category $\cat$ as above, 
one may make a connected category $\cat^+$ by adding a trivial group
to $\cat$ as an initial object (i.e., add one new object, a trivial
group, and one morphism from this object to every other object).  
The natural map from $I(\cat^+)$ to
$I(\cat)$ is an isomorphism, except in degree zero.

The analogue of Corollary~2 in this generality, for which $\cat$ need
not be assumed to be connected, is: 
\proclaim Corollary $2'$.  Let $\cat$ be a finite category of finite
groups and injective homomorphisms.  Then 
$\prod_{Q\in\cat}H^*(Q)$ is a finite module for $I(\cat)$.

\noindent
7) The following instance of Theorem~1 seems worthy of special note.  
Let $P$ be an elementary abelian 2-group of rank $n$, let $\cat$ be
the category whose only object is $P$ and whose morphisms are the 
group $GL(n,\ftwo)$.  Then $H^*(B\Gamma)$ is a ring whose radical 
is invariant under the action of the Steenrod algebra, and 
$H^*(B\Gamma)/\sqrt 0$ is isomorphic to the Dickson algebra 
$D_n= \ftwo[x_1,\ldots,x_n]^{GL(n,\ftwo)}$.  On the other hand it is 
known that for $n\geq 6$, $D_n$ itself cannot be the cohomology of 
any space [\smsw].  

\medskip
\noindent
{\bf Acknowledgements.} This work was done when the authors met at the
Centre de Recerca Matem\`atica.  The authors gratefully acknowledge the
hospitality and support of the CRM.  The first- and second-named
authors were supported by E.C. Leibniz Fellowships, at the
Max-Planck-Institut f\"ur Mathematik and at the CRM respectively.

\bigskip
\goodbreak
\leftline{\bf References.}
\par
\nobreak 
\frenchspacing
\smallfont
\noindent
\item{[\br]} K. S. Brown, Cohomology of Groups, Graduate Texts in
Mathematics {\smallbold 87}, Springer-Verlag, Heidelberg (1982).  

\item{[\ce]} H. Cartan and S. Eilenberg, Homological algebra, 
Princeton Univ. Press, Princeton (1955).

\item{[\wdd]} W. Dicks and M. J. Dunwoody, Groups acting on graphs,
Cambridge Studies in Advanced Mathematics {\smallbold 17}, Cambridge 
Univ. Press, Cambridge (1989).  

\item{[\ep]} L. Evens and S. Priddy, The ring of universally stable
elements, Quart. J. Math. Oxford (2) {\smallbold 40} (1989), 399--407.  

\item{[\ev]} L. Evens, Cohomology of groups, Oxford Mathematical Monographs,
Clarendon Press, New York (1991).  

\item{[\se]} J.-P. Serre, Trees, Springer-Verlag, Heidelberg (1980).  

\item{[\smsw]} L. Smith and R. M. Switzer, Realizability and
non-realizability of the Dickson algebras as cohomology rings, Proc.
Amer. Math. Soc. {\smallbold 89} (1983), 303--313.  

\bigskip

\noindent 
I. J. Leary, Max-Planck-Institut f\"ur Mathematik, 53225 Bonn,
Germany.  

From Jan. 1996: Faculty of Math. Studies, Univ. of Southampton,
Southampton SO17 1BJ, England.  

\noindent 
B. Schuster, CRM, Institut d'Estudis Catalans, E-08193 Bellaterra
(Barcelona), Spain.  

From Jun. 1996: Fachbereich 7 Mathematik, Bergische Universit\"at 
Wuppertal, Gau\ss str. 20, 

42097 Wuppertal, Germany.  

\noindent 
N. Yagita, Faculty of Education, Ibaraki University, Mito, Ibaraki,
Japan.  

\end

Now let $\phi_1,\ldots,\phi_N$ be the
morphisms in the category $\cat$.  Let $Q(i)\leq P$ be the domain of
the homomorphism $\phi_i$.  For each object $Q_j$
of $\cat$, pick $m(j)$ such that $\phi_{m(j)}$ is a morphism from
$Q_j$ to $P$.  To simplify the notation we shall assume that each 
$\phi_{m(j)}$ may be chosen to be the inclusion of $Q_j$ in $P$.  
Let $\Gamma_0$ be the group $P$, and inductively 
define $\Gamma_i$ for $1\leq i\leq N$ by 
$$\Gamma_i=\cases{
\Gamma_{i-1}&if there exists $j$ such that $m(j)=i$,\cr
\Gamma_{i-1}*_{\phi_i,t_i}&otherwise.\cr }$$
In the group $\Gamma_i$, conjugation by $t_i$ as a map from $Q(i)$ is
equal to the map $\phi_i$.  Now let $\Gamma$ be $\Gamma_N$.  

There are many (equivalent) ways to calculate $H^*(\Gamma)$.  Since
the pair of subgroups associated to each of the HNN-extensions made
above were always contained in $P$, they could all be performed
simultaneously, so that $\Gamma$ is the fundamental group of a graph
of groups ([\wdd], p.11) with one vertex $P$, 
and edges the groups $Q(i)$ for those 
$i$ not of the form $i=m(j)$.  This gives rise to a Mayer-Vietoris
type long exact sequence of the form given below.  
$$\cdots\rightarrow
H^k(\Gamma)\rightarrow H^k(P) \rightarrow \bigoplus_{i\neq
m(j)}H^k(Q(i))\rightarrow H^{k+1}(\Gamma)\rightarrow \cdots$$
Here the map from $H^*(P)$ to $H^*(Q(i))$ 
is equal to $\res^P_{Q(i)} - \phi^*\circ\res^P_{\phi(Q(i))}$.  

To see directly that there is a Mayer-Vietoris sequence as above, 
take a model for $B\Gamma$ constructed from a model for $BP$ and 
for each $i$ not of the form $m(j)$ a copy of $BQ(i)\times I$.
The universal cover of this model for $B\Gamma$ consists of lots of 
copies of $EP$, joined together by lots of copies of the
$EQ(i)\times I$'s.  Contracting all the $EQ$'s down to points
leaves a contractible 1-complex (or tree) with a $\Gamma$-action.  
There is one orbit of vertices, of the form $\Gamma/P$, and one orbit
of edges for each $i\neq m(j)$, of orbit type $\Gamma/Q(i)$.  
The equivariant cohomology spectral sequence for this $\Gamma$-complex
has only two non-zero columns and is equivalent to the exact sequence
above.  

Given the above Mayer-Vietoris sequence, claims a) and b) of the
theorem follow easily.  The image of $H^*(\Gamma)$ in $H^*(P)$ is
equal to the intersection of the kernels of the maps 
$\res^P_{Q(i)} - \phi^*\circ\res^P_{\phi(Q(i))}$, which is by
definition $I(P,\cat)$.  The kernel of $\res^\Gamma_P$ consists of
elements in the image of the connecting homomorphism for the exact
sequence, and the product of any two such elements is zero.  (In terms
of the spectral sequence, the only non-zero columns are $E_r^{0,*}$
and $E_r^{1,*}$, and the product of any two elements of $E_r^{1,*}$
lies in $E_r^{2,*}$ so is zero.)  Claim c) is an immediate consequence
of claims a)~and~b), and claim d) of the Theorem follows from
Theorem~5, together with the remarks made in the last paragraph.  
\qed